\input Tex-document.sty

\font\eightmi = cmmi8
\font\boldmathHuge = cmmib10 scaled 1728
\font\boldmathLARGE = cmmib10 scaled 1440

\def\firstpage{139}

\pageno=139

\def\shorttitle{Analyse $p$-adique et Repr\'{e}sentations
Galoisiennes}
\def\shortauthor{J.-M. Fontaine}

\title{\centerline{Analyse {\boldmathHuge p}-adique et} 
\centerline{Repr\'{e}sentations Galoisiennes}}

\author{Jean-Marc Fontaine\footnote{\eightrm *}{\eightrm Institut
Universitaire de France et UMR 8628 du CNRS, Math\'{e}matique,
Universit\'{e} de Paris-Sud, B\^{a}timent 425, 91405 ORSAY Cedex,
France. E-mail: fontaine@math.u-psud.fr}}

\vskip 7mm

\centerline{\boldnormal Abstract}

\vskip 4.5mm

{\narrower \ninepoint \smallskip The notion  of a $p$-adic de Rham
representation of the
absolute Galois group of a $p$-adic field was introduced about
twenty years ago (see e.g. [Fo93]). Three important results for this
theory have been obtained recently: The structure theorem for the {\it
almost  $C_p$-representations}, the theorem {\it weakly admissible
implies admissible} and the theorem {\it de Rham implies potentially
semi-stable}. The proofs of the first two theorems are closely related to
the study of a new kind of analytic groups, the {\it Banach-Colmez spaces}
and the proof of the third  uses deep results on
$p$-{\it adic differential equations on the Robba ring}.
\smallskip

\vskip 4.5mm \noindent {\bf 2000 Mathematics Subject Classification:} 11F80, 11S25, 12H25,  14G22.

\noindent {\bf Keywords and Phrases:} Galois representations, de
Rham representations, Semi-stable representations, $p$-adic Banach
spaces, $p$-adic differential equations.

}

\vskip 10mm
\def\pc#1{\tenrm#1\sevenrm}
\def\pointir{\unskip . --- \ignorespaces}
\def\Medbreak{\vskip-\lastskip\medbreak}

\def\section#1#2{\vfill\eject\ni{\bf {#1} -- {#2}}\bk}

\def\subsection#1#2{\bk\ni{\bf {#1} -- {#2}}\mk}

\long\def\th#1 #2\enonce#3\endth{\Medbreak {\pc#1}
{#2\unskip}\pointir{\it #3}\medskip}

\long\def\th#1 #2\enonce#3\endth{%
\Medbreak {\pc#1} {#2\unskip}\pointir{\it #3}\medskip}

\long\def\thm#1 #2\enonce#3\endthm{%
\Medbreak {\pc#1} {#2\unskip}\pointir{\it #3}\medskip}

\font\tenmsb=msbm10
\font\sevenmsb=msbm7
\font\fivemsb=msbm5
\newfam\msbfam
\textfont\msbfam=\tenmsb
\scriptfont\msbfam=\sevenmsb
\scriptscriptfont\msbfam=\fivemsb

\def\Bbb#1{{\fam\msbfam\relax#1}}

  \let\ni=\noindent
\let\sk=\smallskip \let\mk=\medskip \let\bk=\bigskip
\def\zp {{\Bbb Z}_p} \def\qp{{\Bbb Q}_p} \def\q {{\Bbb Q}} \def\z {\Bbb Z}
\def\bpc{B^+_{cris}}\def\bc{B_{cris}}\def\bdr{B_{dR}}\def\bpdr{B^+_{dR}}
\def\bst{B_{st}} 
\def\fil {{\rm Fil}}
\let\ss=\subset \def\n {\Bbb N} \def\va#1{|#1|}
\def\limproj{\mathop{\oalign{\rm{lim}\cr
\hidewidth$\longleftarrow$\hidewidth\cr}}}
 \def\ext {{\rm Ext}}
\def\bark {\overline K}

  \def\oc {{\cal O}_C}

\def\ac {A_{cris}}
\def\repp {{\rm Rep}_{\qp}(G_K)}
\def\ki {K_{\infty}}
\def\mod {{\rm mod}\ }
\def\bg {{\cal B}(G_K)} \def\rcg {{\rm Rep}_C(G_K)} \def\cg {{\cal C}(G_K)}
\def\barkk {\overline k} \def\no#1{||{#1}||}
\def\spc {{\rm Spm}_C}
\let\f=\rightarrow

\head{1. Repr\'{e}sentations {\boldmathLARGE p}-adiques}

\noindent {\bf 1.1. --}
Dans tout ce qui suit, $K$ est un corps de caract\'{e}ristique $0$,
complet pour une valuation discr\`{e}te, \`{a} corps r\'{e}siduel $k$
parfait de caract\'{e}ristique $p>0$. On choisit une cl\^{o}ture
alg\'{e}brique $\bark$ de $K$, on note $C$ son compl\'{e}t\'{e} et
$\va {\ }_p$ la valeur absolue de $C$ normalis\'{e}e par $\va p _p = p^{-1}$.
On pose $G_K={\rm Gal}(\bark/K)$.

{\it Une repr\'{e}sentation banachique} ({\it de} $G_K$) est un espace de
Banach $p$-adique muni d'une action lin\'{e}aire et continue de
$G_K$. Avec comme morphismes les applications $\qp$-lin\'{e}aires
continues $G_K$-\'{e}quivariantes, ces repr\'{e}sentations forment
une cat\'{e}gorie additive $\qp$-lin\'{e}aire $\bg$.

Une {\it $C$-repr\'{e}sentation} ({\it de} $G_K$) est un $C$-espace
vectoriel de dimension finie muni d'une action semi-lin\'{e}aire et
continue de $G_K$.  Lorsque $k$ est fini, la cat\'{e}gorie $\rcg$ des
$C$-repr\'{e}sentations s'identifie \`{a} une sous-cat\'{e}gorie pleine
de $\bg$~:

\thm PROPOSITION  [Fo00]
\enonce
Supposons $k$ fini. Si $W_1$ et $W_2$ sont des $C$-repr\'{e}senta\-tions,
toute application $\qp$-lin\'{e}aire continue
$G_K$-\'{e}quivariante de $W_1$ dans $W_2$ est $C$-lin\'{e}aire.
\endthm

Disons que deux repr\'{e}sentations banachiques $S_1$ et $S_2$ sont {\it
presque isomorphes} s'il existe un triplet $(V_1,V_2,\alpha)$ o\`{u}
$V_i$ est un sous-$\qp$-espace vectoriel de dimension finie de $S_i$, stable
par $G_K$, et o\`{u} $\alpha:S_1/V_1\f S_2/V_2$ est un
isomorphisme (dans $\bg$). Une {\it presque-$C$-repr\'{e}\-sen\-tation}
({\it de} $G_K$) est une repr\'{e}sentation banachique qui est
presque isomorphe \`{a} une $C$-repr\'{e}sentation. On note $\cg$ la
sous-cat\'{e}gorie pleine de $\bg$ dont les objets sont les
presque-$C$-repr\'{e}sentations. Elle contient la cat\'{e}gorie
$\rcg$ et la cat\'{e}gorie $\repp$ des {\it repr\'{e}sentations
$p$-adiques de dimension finie} ({\it de} $G_K$) comme sous-cat\'{e}gories
pleines.

\thm TH\'{E}OR\`{E}ME A [Fo02]
\enonce
Supposons $k$ fini. La cat\'{e}gorie $\cg$ est ab\'{e}lienne.
Il existe sur les objets de $\cg$ une unique fonction additive $dh:\!
S\mapsto (d(S),h(S))\in\n\times\z$ telle que $dh(W)=(\dim_{C}W,0)$
si $W$ est une
$C$-repr\'{e}\-sen\-ta\-tion et $dh(V)=(0,\dim_{\qp}V)$ si $V$ est de
dimension finie sur $\qp$.

Si $S$ et $T$
sont des objets de $\cg$, les $\qp$-espaces vectoriels
$\ext^{i}_{\cg}(S,T)$ sont de dimension finie et sont nuls pour
$i\not\in\{0,1,2\}$. On a

\centerline{$\sum_{i=0}^{2}(-1)^{i}{\rm dim}_{\qp}{\rm
Ext}^{i}_{\cg}(S,T)=-[K:\qp]h(S)h(T)$.}
\endthm

\ni {\bf 1.2. --} Soit \footnote{$^{(1)}$} {\eightrm Voir [Fo88a]
(resp. [Fo88b]) pour plus de d\'{e}tails sur la construction de
${\hbox{\eightmi B}}_{\hbox{\sixi dR}}$, ${\hbox{\eightmi
B}}_{\hbox{\sixi cris}}$ et ${\hbox{\eightmi B}}_{\hbox{\sixi
st}}$ (resp. sur les repr\'{e}sentations {\eightmi p}-adiques de
de Rham et potentiellement semi-stables).} $FR$ l'ensemble des
suites $(x^{(n)})_{n\in\n}$ d'\'{e}l\'{e}ments de $C$
v\'{e}rifiant $(x^{(n+1)})^{p}=x^{(n)}$ pour tout $n$. Avec les
lois

\centerline{$(x+y)^{(n)}=\lim_{m\mapsto\infty}(x^{(n+m)}+y^{(n+m)})^{p
^{m}}\hbox{
et }(xy)^{(n)}=x^{(n)}y^{(n)}$}

\ni c'est un corps alg\'{e}briquement clos de caract\'{e}ristique $p$,
complet pour la valeur absolue d\'{e}finie par $\va x =\va{x^{(0)}}_p$ et
on note $R$ l'anneau de la valuation. Son corps r\'{e}siduel
s'identifie au corps r\'{e}siduel $\barkk$ de $\bark$. L'anneau $W(R)$
des vecteurs de Witt \`{a} coefficients dans $R$ est int\`{e}gre.
Choisissons $\varepsilon,\pi\in R$ v\'{e}rifiant
$\varepsilon^{(0)}=1$, $\varepsilon^{(1)}\not=1$ et $\pi^{(0)}=p$ et,
pour tout $a\in R$ notons $[a]=(a,0,0,\ldots)$ son repr\'{e}sentant de
Teichm\"{u}ller dans $W(R)$. L'application
$\theta:W(R)\f\oc$ qui envoie $(a_0,a_1,\ldots)$ sur
$\sum_{n\in\n}p^{n}a_n^{(n)}$ est un homomorphisme d'anneaux  surjectif
dont le noyau est l'id\'{e}al principal engendr\'{e} par $\xi=[\pi]-p$. On
note encore $\theta:W(R)[1/p]\f C$ l'application d\'{e}duite
en rendant $p$ inversible. Rappelons que
$\bpdr=\limproj_{n\in\n}W(R)[1/p]/(\xi^{n})$ et que le corps
$\bdr$ des p\'{e}riodes $p$-adiques est le corps des fractions de
$\bpdr$. Toute unit\'{e} $a$ de $R$ s'\'{e}crit d'une mani\`{e}re
unique sous la forme $a=a_0a^{+}$ avec $a_0\in\bar k$ et $\va
{a^{+}-1}<1$, la s\'{e}rie
$\sum_{n=1}^{+\infty}(-1)^{n+1}([a^{+}]-1)^{n}/n$ converge dans
$\bpdr$ vers un \'{e}l\'{e}ment not\'{e} $\log[a]$~; on pose
$t=\log[\varepsilon]$. On a  $\bdr=\bpdr[1/t]$.
On note $\ac$ le s\'{e}par\'{e} compl\'{e}t\'{e} pour la topologie
$p$-adique de la sous-$W(R)$-alg\`{e}bre de $W(R)[1/p]$ engendr\'{e}e
par les $\xi^{m}/m!$ pour $m\in\n$. Alors $\ac$ s'identifie
\`{a} un sous-anneau de $\bpdr$ contenant $t$ et on pose
$\bpc=\ac[1/p]$ et $\bc=\bpc[1/t]\ss\bdr$. La s\'{e}rie
$\sum_{n=1}^{+\infty}(-1)^{n+1}\xi^{n}/np^{n}$ converge dans
$\bpdr$ vers un \'{e}l\'{e}ment $\log[\pi]=\log([\pi]/p)$ et on note
$\bst$ la sous-$\bc$-alg\`{e}bre de $\bdr$ engendr\'{e}e par
$\log[\pi]$. Pour tout $b\in R$ non nul, il existe $r,s\in\z$, avec $s\geq
1$  et une unit\'{e} $a$ de $R$ tels que $b^{s}=\pi^{r}a$ et on pose
$\log[b]=(r\log[\pi]+\log[a])/s$. On a $\bst=\bc[\log[b]]$
d\`{e}s que $b$ n'est pas une unit\'{e}.

Soit ${\cal L}$ l'ensemble des extensions finies de $K$ contenues
dans $\bark$. Pour tout $L\in{\cal L}$, on pose $G_L={\rm Gal}(\bark/L)$
et on note $L_0$ le corps des fractions de l'anneau des vecteurs de
Witt \`{a} coefficients dans le corps r\'{e}siduel de $L$.
Le corps $\bark$ se plonge de fa\c{c}on naturelle dans $\bpdr$ et
l'action de $G_K$ s'\'{e}tend de fa\c{c}on naturelle \`{a} $\bdr$,
l'anneau $\bst$ est stable par $G_K$.  Pour tout $L\in{\cal L}$,
on a $(\bdr)^{G_L}=L$ tandis que $(\bst)^{G_L}=L_0$
et l'application naturelle
$L\otimes_{L_0}\bst\f \bdr$ est injective.

Pour toute repr\'{e}sentation $p$-adique $V$ de
$G_K$ de dimension finie $h$ sur $\qp$, on pose $D_{dR}(V) \allowbreak
=(\bark\otimes_{\qp}V)^{G_K}$,
$D_{st}(V)=(\bst\otimes_{\qp}V)^{G_K}$ et, si $L$ est une extension finie
de $K$ contenue dans $\bark$, $D_{st,L}(V)=(\bst\otimes_{\qp}V)^{G_L}$.
On a $\dim_K D_{dR}(V)\leq h$ et on dit que {\it $V$ est de de Rham}
si on a l'\'{e}galit\'{e}. C'est le cas si $\dim_{K_0}D_{st}(V)=h$
auquel cas on dit que {\it $V$ est semi-stable}. C'est aussi le cas
s'il existe $L\in{\cal L}$ tel que $\dim_{L_0}D_{st,L}(V)=h$, auquel
cas on dit que {\it $V$ est potentiellement semi-stable}, ou si
l'on veut pr\'{e}ciser $L$, que {\it $V$ est $L$-semi-stable}.

\thm TH\'{E}OR\`{E}ME B
\enonce
Toute repr\'{e}sentation $p$-adique de $G_K$ qui est de de Rham est
potentiellement semi-stable.
\endthm
Soit $\bark\bst$ le plus petit sous-anneau de $\bdr$ contenant
$\bark$ et $\bst$. Ce th\'{e}or\`{e}me revient \`{a} dire que, pour toute
repr\'{e}sentation de de Rham $V$, l'inclusion
$((\bark\bst)\otimes_{\qp}V)^{G_K}\ss (\bdr\otimes_{\qp}V)^{G_K}$ est une
\'{e}galit\'{e}. Berger [Be02] en a ramen\'{e} la preuve
\`{a} un r\'{e}sultat sur les \'{e}quations diff\'{e}rentielles
$p$-adiques, r\'{e}sultat qui a ensuite \'{e}t\'{e} prouv\'{e}
ind\'{e}pendamment par Andr\'{e} [An02], Kedlaya [Ke02] et Mebkhout
[Me02], voir \S 3.

L'un des int\'{e}r\^{e}ts de ce th\'{e}or\`{e}me est que l'on dispose
d'une description alg\'{e}bri\-que {\it explicite} de la cat\'{e}gorie
des repr\'{e}sentations potentiellement semi-stables.
Le Frobenius usuel sur $W(R)$ s'\'{e}tend de fa\c{c}on naturelle
en un endomorphisme $\varphi$ de l'anneau $\bst$ (on a $\varphi t=pt$
et $\varphi(\log[\pi])=p\log[\pi]$). Il existe une unique $\bc$-d\'{e}rivation
$N:\bst\f\bst$ telle que $N(\log[\pi])=-1$. L'action de $\varphi$ et de $N$
commutent \`{a} celle de $G_K$ et $N\varphi=p\varphi N$.

Soit $L\in{\cal L}$ telle que $L/K$ est galoisienne. Pour toute
repr\'{e}sentation $p$-adique $V$ de $G_K$, $D_{st,L}(V)$
est un  {\it $(\varphi,N,{\rm Gal}(L/K))$-module filtr\'{e}} de
dimension finie, i.e. c'est un $L_0$-espace vectoriel $D$ de dimension finie,
muni

-- de deux applications $\varphi:D\f D$, $N:D\f D$, la premi\`{e}re
semi-lin\'{e}aire relativement \`{a} la restriction de $\varphi$ \`{a}
$L_0$ et bijective, la deuxi\`{e}me lin\'{e}aire, v\'{e}rifiant
$N\varphi=p\varphi N$,

-- d'une action semi-lin\'{e}aire de ${\rm Gal}(L/K)$, commutant
\`{a} $\varphi$ et \`{a} $N$,

-- d'une filtration index\'{e}e par $\z$, d\'{e}croissante,
exhaustive et s\'{e}par\'{e}e, du $K$-espace vectoriel
$D_K=(L\otimes_{L_0}D)^{{\rm Gal}(L/K)}$ (si $D=D_{st,L}(V)$, on a $D_K\ss
(\bdr\otimes_{\qp}V)^{G_K}$ et, pour tout $i\in\z$,
$\fil^{i}D_K=D_K\cap(\bpdr t^{i}\otimes_{\qp}V)^{G_K}$).

\sk
On pose
$t_H(D)\!=\!\sum_{i\in\z}i.\dim_{K}\fil^{i}D_K/\fil^{i+1}D_K$. Si
$D=\oplus_{\alpha\in\q}D_{\alpha}$  est
la d\'{e}composition suivant les pentes du $\varphi$-isocristal
sous-jacent, on pose aussi
$t_N(D\!)=\!\sum_{\alpha\in\q}\alpha.\dim_{L_0}D_{\alpha}$.  On dit que
$D$ est {\it admissible} si

a) on a $t_H(D)= t_N(D)$,

b) pour tout sous-$L_0$-espace vectoriel $D'$ de $D$, stable par
$\varphi,N$ et ${\rm Gal}(L/K)$, on a $t_H(D')\leq t_N(D')$ (on a muni
$D'_K\ss D_K$ de la filtration induite).

\def\rst {{\rm Rep}_{st,L}(G_K)}
\thm TH\'{E}OR\`{E}ME C [CF00]
\enonce
Soit $L\ss\bark$ une extension finie galoisienne de
$K$.

i) Pour toute repr\'{e}sentation $L$-semi-stable $V$, $D_{st,L}(V)$
est admissible.

ii) Le foncteur qui \`{a} $V$ associe $D_{st,L}(V)$ induit une \'{e}quivalence\footnote{$^{(2)}$} {\eightrm C'est
m\^{e}me une \'{e}quivalence de cat\'{e}gories tannakiennes, cf. [Fo88b].} entre la sous-cat\'{e}gorie pleine
$\rst$ de $\repp$ dont les objets sont les repr\'{e}sentations $L$-semi-stables et la cat\'{e}gorie  des
$(\varphi,N,{\rm Gal}(L/K))$-modules filtr\'{e}s admissibles.
\endthm

\ni {\bf Remarque.} Il \'{e}tait jusqu'\`{a} pr\'{e}sent d'usage
[Fo88b] d'appeler {\it faiblement admissible} ce que nous appelons
ici {\it admissible}. On savait ({\it loc.cit.}, th.5.6.7) que
$D_{st,L}$ induit une \'{e}quivalence entre la cat\'{e}gorie
$\rst$ et une sous-cat\'{e}gorie pleine de la cat\'{e}gorie des
modules filtr\'{e}s (faiblement) admissibles~; on conjecturait que
ce foncteur est essentiellement surjectif et c'est ce qui est
prouv\'{e} dans [CF00].

\head {2. Espaces de Banach-Colmez{\footnote{$^{(3)}$}{\eightrm C'est en cherchant \`{a} prouver le
th\'{e}or\`{e}me C que j'ai \'{e}t\'{e} conduit \`{a} m'int\'{e}resser aux presque {\eightmi
C}-repr\'{e}senta\-tions. C'est Colmez qui a compris que les propri\'{e}t\'{e}s  dont j'avais besoin provenaient
de structures analytiques. Cela nous a permis de prouver le th\'{e}or\`{e}me C. Colmez a ensuite \'{e}tudi\'{e}
plus en d\'{e}tail ces structures analytiques [Co02]. Ce que je raconte ici est une interpr\'{e}tation, dans le
langage de [Fo02], \S 4, de ces travaux de Colmez et devrait \^{e}tre d\'{e}velopp\'{e} dans [FP02].}}}

\ni {\bf 2.1. --} Une  {\it $C$-alg\`{e}bre de Banach} est une
$C$-alg\`{e}bre norm\'{e}e compl\`{e}te $A$~; son {\it spectre
maximal} est l'ensemble $\spc A$ des sections continues  $s : A\f
C$ du morphisme structural. Si $f\in A$ et $s\in\spc A$, on pose
$f(s)=s(f)$. Une {\it $C$-alg\`{e}bre spectrale} est une
$C$-alg\`{e}bre de Banach $A$ telle que la norme est la norme
spectrale, i.e. telle que, pour tout $f\in A$, $\no f
=sup_{s\in\spc A}\va {f(s)}_p$~; dans ce cas,  $\spc A$ est un
espace m\'{e}trique complet (la distance \'{e}tant d\'{e}finie par
$d(s_1,s_2)=\sup_{\no f\leq 1}\va {f(s_1)-f(s_2)} _p$). Avec comme
morphismes les homomorphismes continus de $C$-alg\`{e}bres, les
$C$-alg\`{e}bres spectrales forment une cat\'{e}gorie. La {\it
cat\'{e}gorie des vari\'{e}t\'{e}s spectrales affines sur $C$} est
la cat\'{e}gorie oppos\'{e}e.

Un {\it groupe spectral commutatif affine sur $C$} est un objet en
groupes
commutatifs dans la cat\'{e}gorie des vari\'{e}t\'{e}s spectrales
affines sur $C$. Ces groupes forment, de fa\c{c}on \'{e}vidente, une
cat\'{e}gorie additive qui a des limites projectives finies. Le
foncteur qui \`{a} un groupe spectral commutatif affine associe le
groupe topologique sous-jacent est fid\`{e}le. Si ${\cal S}=\spc A$
est un groupe
spectral commutatif affine, {\it un sous-groupe spectral} du groupe
topologique sous-jacent est un sous-groupe ${\cal T}$ qui admet une
structure de groupe spectral (n\'{e}cessairement unique) telle que
l'inclusion ${\cal T}\f {\cal S}$ est un morphisme de groupes spectraux.

Soit $S$ un espace de Banach ($p$-adique) et ${\cal S}_0$ la
boule unit\'{e}. Un {\it r\'{e}seau} de $S$ est un
sous-$\zp$-module ${\cal S}$ qui est tel que l'on peut trouver
$r,s\in\z$ v\'{e}rifiant $p^{r}{\cal S}_0\ss {\cal S}\ss p^{s}{\cal
S}_0$. Il revient au m\^{e}me de dire qu'il existe une norme
\'{e}quivalente \`{a} la norme donn\'{e}e pour laquelle ${\cal S}$ est
la boule unit\'{e}.

Une {\it $C$-structure
analytique sur $S$} est la donn\'{e}e d'un $C$-groupe spectral
commutatif affine ${\cal S}$ et d'un homomorphisme continu du groupe
topologique sous-jacent \`{a} ${\cal S}$ dans $S$ dont l'image est
un r\'{e}seau et le noyau un $\zp$-module de type fini.
On dit que deux $C$-structures
analytiques ${\cal S}$ et ${\cal T}$ sur $S$ sont {\it
\'{e}quivalentes} si ${\cal S}\times_{S}{\cal T}$ est un
sous-groupe spectral de ${\cal S}\times{\cal T}$. Un {\it (espace
de) Banach analytique (sur $C$)} est la donn\'{e}e d'un espace de
Banach muni d'une classe d'\'{e}quivalence de $C$-structures
analytiques (on les appelle les structures {\it admissibles} de $S$).
On dit que $S$ est {\it effectif} s'il existe une structure
admissible ${\cal S}$ telle que l'application ${\cal S}\f S$ est
injective.

Un {\it morphisme de Banach analytiques} $f : S\f T$ est une application
$\qp$-lin\'{e}aire continue telle qu'il existe des structures
admissibles ${\cal S}$ de $S$ et ${\cal T}$ de $T$ et un
morphisme ${\cal S}\f {\cal T}$ qui induit $f$. Les Banach
analytiques forment une cat\'{e}gorie additive ${\cal BA}_C$.

Si $S$ est un Banach analytique et si $V$ est un sous-$\qp$-espace
vectoriel de dimension finie, le quotient $S/V$ a une structure
naturelle  de Banach analytique. On dit que deux Banach analytiques
$S_1$ et $S_2$ sont {\it presque isomorphes} s'il existe des
sous-$\qp$-espaces vectoriels de dimension finie $V_1$ de $S_1$ et
$V_2$ de $S_2$ et un isomorphisme $S_1/V_1\f S_2/V_2$ (de Banach
analytiques).

Le groupe sous-jacent \`{a} $\oc$ a une structure naturelle de groupe
spectral commutatif affine~: on a $\oc=\spc C\{X\}$ o\`{u} $C\{X\}$
est l'alg\`{e}bre de Tate des s\'{e}ries formelles \`{a} coefficients dans $C$
en l'ind\'{e}termin\'{e}e $X$ dont le terme g\'{e}n\'{e}ral tend vers
$0$. Ceci fait de $C$ un espace de Banach analytique effectif. Un
{\it Banach analytique vectoriel} est un Banach analytique isomorphe \`{a}
$C^{d}$ pour un entier $d$ convenable. Un {\it espace de Banach-Colmez}
est un {\it Banach analytique presque vectoriel}, i.e. un Banach
analytique qui est presque isomorphe \`{a} un Banach analytique
vectoriel. On note ${\cal BC}_C$ la sous-cat\'{e}gorie pleine de ${\cal
BA}_C$ dont les objets sont les espaces de Banach-Colmez.

\thm PROPOSITION (th\'{e}or\`{e}me de Colmez {\footnote{$^{(4)}$}{\eightrm C'est \`{a} peu pr\`{e}s  le
r\'{e}sultat principal de [Co02]. La d\'{e}finition donn\'{e}e par Colmez de ce qu'il appelle les Espaces de
Banach de dimension finie (avec un E majuscule) est l\'{e}g\`{e}rement diff\'{e}rente. Il n'est pas tr\`{e}s
difficile de construire une \'{e}quivalence entre sa cat\'{e}gorie et la n\^{o}tre [FP02].}}) \enonce La
cat\'{e}gorie ${\cal BC}_C$ est ab\'{e}lienne et le foncteur d'oubli de ${\cal BC}_C$ dans la cat\'{e}gorie des
$\qp$-espaces vectoriels est exact et fid\`{e}le. Il existe sur les objets de ${\cal BC}_C$, une unique fonction
additive $dh: S\mapsto (d(S),h(S))\in\n\times\z$ telle que $dh(C^{d})=(d,0)$ et $dh(V)=(0,\dim_{\qp}V)$ si $V$ est
de dimension finie sur $\qp$.
\endthm

\ni {\bf 2.2. --} La meilleure fa\c{c}on de comprendre les
th\'{e}or\`{e}mes A et C c'est d'utiliser le r\'{e}sultat
pr\'{e}c\'{e}dent pour les prouver {\footnote{$^{(5)}$}{\eightrm
C'est ainsi que Colmez red\'{e}montre le th\'{e}or\`{e}me C dans
[Co02]. Moyennant une preuve un peu plus compliqu\'{e}e, on peut
n'utiliser qu'un r\'{e}sultat d'analyticit\'{e} apparemment moins
fort~; c'est ce qu'on fait  pour prouver le th\'{e}or\`{e}me C
dans [CF00] et le th\'{e}or\`{e}me A dans [Fo02].}. Lorsque $k$
est fini, toute presque-$C$-repr\'{e}sen\-ta\-tion est munie
d'unee structure naturelle  d'espace de Banach-Colmez~; toute
application $\qp$-lin\'{e}aire continue $G_K$-\'{e}quivariante
d'une presque $C$-repr\'{e}sentation dans une autre est
analytique. Le fait que $\cg$ est ab\'{e}lienne et l'existence de
la fonction $dh$ r\'{e}sultent alors du th\'{e}or\`{e}me de
Colmez.

Le principe de la preuve du th\'{e}or\`{e}me C est le suivant~: On se
ram\`{e}ne facilement au cas semi-stable, i.e. au cas o\`{u} $L=K$.
Il s'agit de v\'{e}rifier que, si $D$ est un $(\varphi,N)$-module
filtr\'{e} (faiblement) admissible de dimension $h$,
il existe une repr\'{e}sentation $p$-adique $V$ de dimension $h$
telle que $D_{st,K}(V)$ soit isomorphe \`{a} $D$. Une torsion \`{a} la
Tate permet de supposer que $\fil^{0}D_K=D_K$. Notons
$V_{st}^{+,0}(D)$ le $\qp$-espace vectoriel des applications
$K_0$-lin\'{e}aires de $D$ dans $B_{st}\cap\bpdr$ qui commutent \`{a}
l'action de $\varphi$ et de $N$ et $V_{st}^{+,1}$ le quotient du $K$-espace
vectoriel des applications $K$-lin\'{e}aires de $D_K$ dans $\bpdr$ par
le sous-espace des applications qui sont compatibles avec la filtration.
On commence par
v\'{e}rifier que le noyau $V_{st}^{\ast}(D)$ de l'application \'{e}vidente
$ \beta : V_{st}^{+,0}(D)\f V_{st}^{+,1}(D)$ est un $\qp$-espace
vectoriel de dimension finie $\leq h$ et que, s'il est de dimension
$h$, alors la repr\'{e}sentation duale $V_{st}(D)$ est semi-stable et
$D$ est isomorphe \`{a} $D_{st}(V_{st}(D))$. La th\'{e}orie des espaces
de Banach analytiques permet de munir $V_{st}^{+,0}(D)$ et
$V_{st}^{+,1}(D)$ d'une structure d'espace de Banach-Colmez et on a
$dh(V_{st}^{+,0}(D))=(t_N(D),h)$ tandis que
$dh(V_{st}^{+,1}(D))=(t_H(D),0)$. Il suffit alors de v\'{e}rifier que
l'application $\beta$ est analytique. Comme $t_H(D)=t_N(D)$,
l'additivit\'{e} de $dh$ implique
que $\beta$  est surjective et que
$dh(V_{st}^{\ast}(D))=(0,h)$, ce qui signifie bien  que
$\dim_{\qp}V_{st}^{\ast}(D)=h$.

\head {3. Equations
diff\'{e}rentielles}

\ni {\bf 3.1. --} Soit $A$ un anneau commutatif et $d:A\f\Omega_A$ une d\'{e}rivation de $A$ dans un $A$-module
$\Omega_A$. Ici, un $A$-module \`{a} connexion (sous-entendu relativement \`{a} $d$) est un $A$-module libre de
rang fini ${\cal D}$ muni d'une application $\nabla : {\cal D}\f {\cal D}\otimes\Omega_A$ v\'{e}rifiant la
r\`{e}gle de Leibniz. On dit que ce module est {\it trivial} s'il est engendr\'{e} par le sous-groupe ${\cal
D}_{\nabla=0}$ des {\it sections horizontales}.

Pour tout corps $L$ de caract\'{e}ristique $0$,
complet pour une valuation discr\`{e}te, notons (cf. par
exemple [Ts98], \S 2) ${\cal R}_{x,L}$ l'{\it anneau de Robba de}
$L$ (ou `'anneau des fonctions
analytiques sur une couronne d'\'{e}paisseur nulle''),
c'est-\`{a}-dire l'anneau
des s\'{e}ries $\sum_{n\in\z}a_nx^{n}$ \`{a} coefficients dans $L$
v\'{e}rifiant

\centerline{$\forall s<1, \va{a_n}s^{n}\mapsto 0\hbox{ si }n\mapsto +\infty
\hbox{ et }\exists r<1\hbox{ tel que }\va{a_n}r^{n}\mapsto 0\hbox{ si }
n\mapsto -\infty$.}

Le sous-anneau ${\cal E}_{x,L}^{\dag}$ de ${\cal R}_{x,L}$ des fonctions
$\sum a_nx^{n}$ telles que les $a_n$ sont born\'{e}s est un corps muni
d'une valuation discr\`{e}te (d\'{e}finie par
$\va{\sum a_nx^{n}}={\rm sup}\va{a_n}$)
qui n'est pas complet mais est hens\'{e}lien. Son
corps r\'{e}siduel s'identifie au corps des s\'{e}ries formelles
$E=k_L((x))$ o\`{u} $k_L$ d\'{e}signe le corps r\'{e}siduel de $L$.
Pour toute extension finie s\'{e}parable $F$ de $E$, il existe
une, unique \`{a} isomorphisme unique pr\`{e}s, extension non
ramifi\'{e}e ${\cal E}_F^{\dag}$ de ${\cal E}_{x,L}^{\dag}$ de
corps r\'{e}siduel $F$ . Posons ${\cal
R}_{F}={\cal R}_{x,L}\otimes_{{\cal E}^{\dag}_{x,L}}{\cal
E}_F^{\dag}$.
Si $k_F$ d\'{e}signe le corps r\'{e}siduel de $F$, $L'$ l'unique
extension non ramifi\'{e}e de $L$ de corps r\'{e}siduel $k_F$ et si
$x'$ est un rel\`{e}vement dans l'anneau des entiers de
${\cal E}^{\dag}_F$ d'une uniformisante de $F$, l'anneau ${\cal R}_F$
s'identifie \`{a} l'anneau de Robba ${\cal R}_{x',L'}$.

Notons $\Omega^{1}_{{\cal R}_{x,L}}$ le ${\cal R}_{x,L}$-module libre
de rang $1$ de base $dx$, solution du probl\`{e}me universel pour les
d\'{e}rivations continues en un sens \'{e}vident. Les modules \`{a}
connexion sur l'anneau ${\cal R}_{x,L}$} forment une cat\'{e}gorie artinienne.
Si ${\cal D}$ est un objet de cette cat\'{e}gorie, on dit qu'il est
{\it unipotent}
si son semi-simplifi\'{e} est trivial. On dit qu'il est
{\it quasi-unipotent} s'il existe  une extension finie s\'{e}parable
$F$ de $k((X))$ telle que le module \`{a} connexion sur ${\cal R}_F$
d\'{e}duit de ${\cal D}$ par extension des scalaires soit unipotent.

Pour tout $z$ dans l'aneeau des entiers de  ${\cal E}_{x,K_0}^{\dag}$,
il existe un unique endomorphisme continu $\varphi$ de ${\cal R}_{x,K_0}$
qui prolonge le Frobenius absolu sur $K_0$ et v\'{e}rifie
$\varphi(x)=x^{p}+pz$~; on appelle {\it Frobenius} un tel
endomorphisme. Pour un tel $\varphi$, on note encore $\varphi:
\Omega^{1}_{{\cal R}_{x,K_0}}\f \Omega^{1}_{{\cal R}_{x,K_0}}$ l'application
induite.
Soit ${\cal D}$ un module \`{a} connexion sur ${\cal R}_{x,K_0}$. Une
{\it structure de Frobenius sur  ${\cal D}$} consiste en la
donn\'{e}e d'un Frobenius $\varphi$ sur ${\cal R}_{x,K_0}$ et d'une
application $\varphi$-semi-lin\'{e}aire $\varphi_{{\cal D}}:{\cal D}\f{\cal D}$
commutant \`{a} $\nabla$.

\thm TH\'{E}OR\`{E}ME (Andr\'{e}, Kedlaya, Mebkhout \footnote
{$^{(6)}$}
{\eightrm Crew [Cr98] a sugg\'{e}r\'{e} que ce th\'{e}or\`{e}me pouvait
\^{e}tre vrai~; il a \'{e}t\'{e}
prouv\'{e} ind\'{e}pen\-damment par Andr\'{e} [An02], Mebkhout [Me02]
et Kedlaya [Ke01]. Pour Andr\'{e} comme pour Mebkhout, c'est un cas
particulier d'un r\'{e}sultat plus g\'{e}n\'{e}ral dont la preuve
repose sur la th\'{e}orie de Christol-Mebkhout [CM]. La preuve de
Kedlaya est plus directe~: elle utilise une classification \`{a} la
Dieudonn\'{e}-Manin des modules munis d'un Frobenius
pour se ramener \`{a} un r\'{e}sultat de Tsuzuki [Ts98]. Voir [Co01] pour une
\'{e}tude comparative  plus d\'{e}taill\'{e}e.})
\enonce
Tout module \`{a} connexion sur ${\cal R}_{x,K_0}$ qui admet
une structure de Frobenius est quasi-unipotent.
\endthm
  Avant de montrer comment Berger [Be2] d\'{e}duit le
  th\'{e}or\`{e}me B de cet \'{e}nonc\'{e}, rappelons quelques
  r\'{e}sultats de [Fo00], [Fo90] et [CC98] (cf. aussi [Co98]). Dans
tout ce qui suit,
  $V$ est une repr\'{e}sentation $p$-adique de $G_K$ de dimension
  finie $h$.
\sk

\def\oe#1{{\cal O}_{{\cal E}_{#1}}}
\def\oenr {{\cal O}_{\widehat{\cal E}^{nr}}}
\def\enr {\widehat{\cal E}^{nr}}

\ni {\bf 3.2. --} Soit $\ki$ le sous-corps de $\bark$ engendr\'{e} sur $K$ par les racines de l'unit\'{e} d'ordre
une puissance de $p$. Posons $H_K={\rm Gal}(\bark/\ki)$ et $\Gamma_K=G_K/H_K$. En utilisant la th\'{e}orie de Sen
[Se80], on montre [Fo00] que l'union $\Delta_{dR}(V)$ des sous-$\ki[[t]]$-modules de type fini de
$(\bdr\otimes_{\qp}V)^{H_K}$ stables par $\Gamma_K$ est un  $\ki((t))$-espace vectoriel de dimension $h$ et qu'il
existe une unique connexion
$$\nabla : \Delta_{dR}(V)\f \Delta_{dR}(V)\otimes dt/t$$
qui a la propri\'{e}t\'{e} que, pour tout sous-$\ki[[t]]$-module de type fini
$Y$ stable par $G_K$, tout entier $r\geq 0$ et tout $y\in Y$, il existe un
sous-groupe ouvert $\Gamma_{r,y}$ de $\Gamma$ tel que, si
$\nabla(y)=\nabla_0(y)\otimes dt/t$, alors

\centerline{$\gamma(y)\equiv \exp(\log\chi(\gamma).\nabla_0)(y)\ \
(\mod t^{r}Y)$,
pour tout $\gamma\in\Gamma_{r,y}$.}

Cette connexion est {\it r\'{e}guli\`{e}re}~: le $\ki[[t]]$-module
$\Delta^{+}_{dR}(V)=(\bpdr\otimes V)\cap\Delta_{dR}(V)$ est un
r\'{e}seau de $\Delta_{dR}(V)$ v\'{e}rifiant
$\nabla(\Delta^{+}_{dR}(V))\ss \Delta^{+}_{dR}(V)\otimes dt/t$.
On a $D_{dR}(V)=(\Delta_{dR}(V))^{\Gamma_K}$. L'action de $\Gamma_K$
est discr\`{e}te sur $\Delta_{dR}(V)_{\nabla=0}$~; on en d\'{e}duit
que $\Delta_{dR}(V)_{\nabla=0}=\ki\otimes_K D_{dR}(V)$ donc que $V$
est de de Rham si et seulement si le module \`{a} connexion
$\Delta_{dR}(V)$ est trivial. Ceci se produit si et seulement s'il
existe un r\'{e}seau (n\'{e}cessairement unique) $\Delta^{0}_{dR}(V)$
de $\Delta_{dR}(V)$
v\'{e}rifiant $\nabla(\Delta^{0}_{dR}(V))\ss \Delta^{0}_{dR}(V)\otimes
dt$.
\sk

\ni {\bf 3.3. --} Rappelons bri\`{e}vement la th\'{e}orie des $(\varphi,\Gamma)$-modules [Fo80]. Soit $\oe 0$
l'adh\'{e}\-rence dans $W(FR)$ de la sous-$W(k)$-alg\`{e}bre engendr\'{e}e par $[\varepsilon]$ et
$1/([\varepsilon]-1)$. C'est un anneau de valuation discr\`{e}te complet dont l'id\'{e}al maximal est engendr\'{e}
par $p$ et dont le corps r\'{e}siduel $E_0$ est le corps des s\'{e}ries formelles $k((\varepsilon-1))$ vu comme
sous-corps ferm\'{e} de $FR$.  Notons $\oenr$ le s\'{e}par\'{e} compl\'{e}t\'{e} pour la topologie $p$-adique de
l'union de toutes les sous-$\oe 0$-alg\`{e}bres finies \'{e}tales de $\oe 0$ contenues dans $W(FR)$. C'est un
anneau de valuation discr\`{e}te complet dont le corps r\'{e}siduel  est une cl\^{o}ture s\'{e}parable $E^{s}$ de
$E_0$. Son corps des fractions $\enr$ s'identifie \`{a} un sous-corps ferm\'{e} du corps  ${\widetilde B}=
W(FR)[1/p]$, stable par l'action de $G_K$ et par le Frobenius $\varphi$. Le corps ${\cal E}_K=(\enr)^{H_K}$ est
une extension finie non ramifi\'{e}e du corps des fractions de $\oe 0$. Son corps r\'{e}siduel $E_K$ est une
extension finie s\'{e}parable de $E_0$~; le corps r\'{e}siduel  $k'$ de $E_K$ est celui de $\ki$.

Alors, $D(V)=(\enr\otimes_{\qp}V)^{H_K}$ est un
{\it $(\varphi,\Gamma_K)$-module sur ${\cal E}_K$}, i.e. un
${\cal E}_K$-espace vectoriel de dimension finie $D$ muni d'un
Frobenius $\varphi$-semi-lin\'{e}aire (que l'on note encore $\varphi$)
et d'une action
semi-lin\'{e}aire continue de $\Gamma_K$ commutant \`{a} l'action de
$\varphi$~; ce $(\varphi,\Gamma_K)$-module est {\it \'{e}tale}, i.e.
il existe un $\oe {K}$-r\'{e}seau ${\cal D}$ de $D$ tel que ${\cal D}$
est le $\oe {K}$-module engendr\'{e} par $\varphi({\cal D})$. La
correspondance $V\mapsto D(V)$ d\'{e}finit une \'{e}quivalence entre
$\repp$ et la cat\'{e}gorie des $(\varphi,\Gamma_K)$-modules \'{e}tales.
\def\enrd {\widehat{\cal E}^{nr,\dag}} \def\edk#1{{\cal E}^{\dag}_{K,#1}}
\def\ekr {{\cal E}_K^{rig}}
\sk \ni {\bf 3.4. --} Il n'y a pas de fl\`{e}che naturelle de $\enr$ dans $\bpdr$, ce qui fait que la comparaison
entre  $\Delta_{dR}(V)$ et $D(V)$ n'est pas si facile. Toutefois, si $a\in R$ est non nul, $[a]\in W(R)\ss\bpdr$
est inversible dans $\bpdr$, ce qui permet de voir $1/[a]=[1/a]$ comme un \'{e}l\'{e}ment de $\bpdr$. Tout
\'{e}l\'{e}ment de ${\widetilde B}$ s'\'{e}crit d'une mani\`{e}re et d'une seule sous la forme
$\sum_{n\gg-\infty}p^{n}[a_n]$, avec les $a_n\in FR$~; notons ${\widetilde B}^{+}_{dR}$ le sous-anneau de
${\widetilde B}$ form\'{e} des s\'{e}ries de ce type qui convergent dans $\bpdr$. L'application ${\widetilde
B}^{+}_{dR}\f\bpdr$ est injective et permet d'identifier ${\widetilde B}^{+}_{dR}$ \`{a} une
sous-$W(R)[1/p]$-alg\`{e}bre de $\bpdr$. Pour tout $r\in\n$, posons $\enrd_r=\enr\cap\varphi^{r}({\widetilde
B}^{+}_{dR})$ et, pour tout $b\in \enrd_r$, notons $\varphi_r(b)$ l'unique $c\in {\widetilde B}^{+}_{dR}\ss\bpdr$
tel que $\varphi^{r}(c)=b$. On a $\enrd_{r}\ss{\enrd}_{r+1}$ ; soit $\enrd$ l'union des $\enrd_r$. Alors ${\cal
E}_{K}^{\dag}=(\enrd)^{H_K}$ est un sous-corps dense de ${\cal E}_K$ stable par $\varphi$. On pose
$D^{\dag}(V)=(\enrd\otimes_{\qp}V)^{H_K}$. On peut le calculer \`{a} partir de $D(V)$~: c'est  l'union des
sous-${\cal E}_{K}^{\dag}$-espaces vectoriels de dimension finie de $D(V)$ stables par $\varphi$. Le r\'{e}sultat
principal de [CC98] est que $V$ est surconvergente, c'est-\`{a}-dire que l'application naturelle ${\cal
E}_K\otimes D^{\dag}(V)\f D(V)$ est un isomorphisme.

Pour tout $r\in\n$, soit $\edk r = (\enrd_r)^{H_K}$. Alors
$D^{\dag}_r(V)=(\enrd_r\otimes_{\qp}V)^{H_K}$ est aussi le plus
grand sous-$\edk r$-module $M$ de type fini de $D(V)$ tel que $\varphi(M)\ss
\edk{r+1}M$. Pour $r$ assez grand, l'application naturelle
${\cal E}_{K}^{\dag}\otimes_{\edk r}D^{\dag}_r(V)\f D^{\dag}(V)$ est un
isomorphisme. Lorsqu'il en est ainsi, on a $\varphi_r(\edk r)\ss
\ki[[t]]$ et

\centerline{$\Delta_{dR}(V)=\ki((t))\otimes_{\edk r}D^{\dag}_r(V) \hbox{ et
donc }
D_{dR}(V)= (\ki((t))\otimes_{\edk r}D^{\dag}_r(V))^{\Gamma_K}.$}

\ni {\bf 3.5. --} Wach [Wa96] a montr\'{e} comment calculer
$D_{st}(V)$ \`{a} partir de $D^{\dag}(V)$ lorsque $V$ est {\it de
hauteur finie}. C'est Berger [Be02] qui a compris comment traiter
le cas g\'{e}n\'{e}ral~: Choisissons un rel\`{e}vement $x$ dans
l'anneau des entiers de ${\cal E}_{K}^{\dag}$ d'une uniformisante
de $E_K$.  Si $K'_0$ d\'{e}signe le corps des fractions de
$W(k')$, ${\cal E}_{K}^{\dag}$ s'identifie pr\'{e}cis\'{e}ment au
sous-anneau ${\cal E}_{x,K'_0}^{\dag}$ de l'anneau de Robba ${\cal
R}_{x,K'_0}$. Ce dernier ne d\'{e}pend pas du choix de $x$ et
nous le notons $\ekr$~; il contient $t= \log[\varepsilon]$. Si $L$
est une extension finie de $K$ contenue dans $\bark$, le corps
$F=(E^{s})^{{\rm Gal}(\bark/L\ki)}$ est une extension finie de
$E_K$ et l'anneau ${\cal E}_L^{rig}$ s'identifie \`{a}  l'anneau
not\'{e} ${\cal R}_F$ au \S 3.1.

Rappelons (\S1.2) que si $u=\log[\varepsilon-1]$,
on a $\bst=\bc[u]$. Les actions de $\varphi$ et de $\Gamma_K$
s'\'{e}tendent de fa\c{c}on \'{e}vidente \`{a} $\ekr$,  \`{a} $\ekr[1/t]$
et \`{a} l'anneau $\ekr[1/t][u]$
des polyn\^{o}mes en $u$ \`{a} coefficients dans $\ekr[1/t]$. Berger
montre que $$D_{cris}(V)=(\ekr[1/t]\otimes_{{\cal
E}_{K}^{\dag}}D^{\dag}(V))^{\Gamma_K}
\hbox{ et }D_{st}(V)=(\ekr[1/t][u]\otimes_{{\cal
E}_K^{\dag}}D^{\dag}(V))^{\Gamma_K}$$
(l'action de $N$ sur $D_{st}(V)$ est la restriction de
$-d/du\otimes id_{D{\dag}(V)}$).
\mk

\ni {\bf 3.6. --} Posons $D=D_K^{rig}(V)=\ekr[1/t]\otimes_{{\cal
E}_{K}^{\dag}}D^{\dag}(V)$. En utilisant l'action de $\Gamma_K$
comme au \S 3.4, on d\'{e}finit une connexion $\nabla : D \f D
\otimes dt/t$ qui commute \`{a} l'action de $\varphi$. Cette
connexion est {\it r\'{e}guli\`{e}re} au sens qu'il existe un
sous-$\ekr$-module $D^{+}$ de $D$, libre de rang $h$, stable par
$\varphi$ et v\'{e}rifiant $\nabla(D^{+})\ss D^{+}\otimes dt/t$
(prendre $D^{+}=\ekr\otimes_{{\cal E}_K^{\dag}}D^{\dag}(V)$). On
v\'{e}rifie que le $\ekr$-module libre $\Omega^{1}_{\ekr}$ admet
$d[\varepsilon]$ comme base. Mais $dt/t=
[\varepsilon]^{-1}/td[\varepsilon]$ et $t$ n'est pas inversible
dans $\ekr$. On d\'{e}duit alors facilement du th\'{e}or\`{e}me
d'Andr\'{e}-Kedlaya-Mebkhout que $V$ est potentiellement
semi-stable si et seulement s'il existe un sous-$\ekr$-module
libre $D^{0}$ de $D$, libre de rang $h$, stable par $\varphi$ et
v\'{e}rifiant $\nabla(D^{0})\ss D^{0}\otimes dt$.

Il ne reste plus qu'\`{a} construire un tel $D^{0}$ lorsque $V$ est de
de Rham. Fixons un entier $r_0\geq 1$ suffisamment grand pour que
$D^{\dag}_{r_0}(V)$ contienne une base de $D^{\dag}(V)$ sur ${\cal
E}_K^{\dag}$ et pour que  $x\in{\cal E}^{\dag}_{r_0}$. Pour tout
$r\geq r_0$ le sous-anneau ${\cal E}_{K,r}^{rig}$ de $\ekr={\cal
R}_{x,K'_0}$ form\'{e} des
$\sum_{n\in\z}a_nx^{n}$ v\'{e}rifiant

\centerline{$\forall s<1, \va{a_n}s^{n}\mapsto 0\hbox{ si }n\mapsto +\infty
\hbox{ et }a_n(\varepsilon^{r}-1)^{n}\mapsto 0\hbox{ si }
n\mapsto -\infty$}

\ni est stable par $\Gamma_K$ et contient $\edk r$. Si
$D_r=D_{K,r}^{rig}(V)={\cal
E}_{K,r}^{rig}\otimes_{\edk r}D_r^{\dag}(V)$, alors $D$ est
la r\'{e}union croissante des $D_r$ et $\varphi(D_r)\ss D_{r+1}$. L'application
$\varphi_r$ induit un homomorphisme de ${\cal E}_{K,r}^{rig}$ dans
$\ki((t))$ et un isomorphisme de
$\ki((t))\otimes_{{\cal E}_{K,r}^{rig}}D_r$
sur $\Delta_{dR}(V)$. L'application $\Phi_r:D_r\f \Delta_{dR}(V)$
qui envoie $a$ sur $1\otimes a$ est injective.

Soit $\Delta_{dR}^{0}(V)$ le sous-$\ki((t))$-module de $\Delta_{dR}(V)$
engendr\'{e} par les sections horizontales. Pour tout $r\geq r_0$,
soit $D_r^{0}=\{a\in D_r\mid \Phi_s(a)\in
\Delta^{0}_{dR}(V)\hbox{ pour tout }s\geq r\}$. On a $D_r^{0}\ss
D_{r+1}^{0}$ et $D^{0}=\cup_{r\geq r_0}D_r^{0}$ est un
sous-$\ekr$-module de $D$, stable par $\varphi$ et v\'{e}rifiant
$\nabla(D^{0})\ss D^{0}\otimes dt$.  Si $V$ est de de Rham on
d\'{e}duit du fait que
$\Delta_{dR}^{0}(V)$ est un r\'{e}seau de $\Delta_{dR}(V)$ que
$D^{0}$ est libre de rang $h$ sur $\ekr$. D'o\`{u} le th\'{e}or\`{e}me B.
\vfill\eject

\parindent 11mm

\head{Bibliographie}

\item{[An02]} Y. Andr\'{e}, {\it Filtration de type Hasse-Arf et monodromie
$p$-adique}, Inv. Math. {\bf 148} (2002), 285--317.

\item{ [Be02]} L. Berger, {\it Repr\'{e}sentations $p$-adiques et
\'{e}quations diff\'{e}rentielles}, Inv. Math. {\bf 148} (2002),
219-284.

\item{ [CM]} G. Christol et Z. Mebkhout, {\it Sur le th\'{e}or\`{e}me de
l'indice des \'{e}quations diff\'{e}ren\-tielles $p$-adiques} I, Ann. Inst.
Fourier {\bf 43} (1993), 1545--1574~; II, Ann. of Maths. {\bf 146}
(1997), 345--410~; III, Ann. of Maths. {\bf 151} (2000), 385--457~;
IV, Inv. Math. {\bf 143} (2001), 629--671.

  \item{ [Co98]} P. Colmez {\it Repr\'{e}sentations $p$-adiques d'un
corps  local}
  in Proceedings of the I.C.M. Berlin,
  vol. II, Documenta Mathematica (1998), 153--162.

  \item{ [Co01]} P. Colmez, {\it Les conjectures de monodromie $p$-adique},
   S\'em. Bourbaki,  exp. 897, novembre 2001.

\item{ [Co02]} P. Colmez, {\it Espaces de Banach de dimension finie},
J. Inst. Math. Jussieu, \`{a} para\^{\i}tre.

\item{ [CC98]} F. Cherbonnier et P. Colmez, {\it Repr\'{e}sentations
$p$-adiques surconvergentes}, Inv. Math. {\bf 133} (1998), 581--611.

\item{ [CF00]} P. Colmez et J.-M. Fontaine, {\it Construction des
repr\'{e}sen\-tations semi-stables}, Inv. Math. {\bf 140} (2000), 1--43.

\item{ [Cr98]} R. Crew, {\it Finiteness theorems for the cohomology of an
overconvergent isocrystal on a curve}, Ann. scient. E.N.S. {\bf 31}
(1998), 717--763.

\item{ [Fo83]} J.-M. Fontaine, {\it Repr\'{e}\-sen\-ta\-tions $p$-adiques},
in Proceedings of the I.C.M.,
War\-sza\-wa, vol. I, Elsevier, Amsterdam (1984), 475--486.

\item{ [Fo88a]} J.-M. Fontaine, {\it Le corps des p\'{e}riodes $p$-adiques},
avec un appendice par Pierre
Colmez, in  P\'{e}riodes $p$-adiques, Ast\'e\-risque {\bf 223},
S.M.F., Paris (1994), 59--111.

\item{ [Fo88b]} J.-M. Fontaine, {\it Repr\'{e}\-sen\-ta\-tions $p$-adiques
semi-stables},
in  P\'eriodes $p$-adiques, Ast\'e\-risque {\bf 223},
S.M.F., Paris (1994), 113--184.

\item{ [Fo90]} J.-M. Fontaine, {\it Repr\'{e}\-sen\-ta\-tions $p$-adiques des
corps locaux}, in the Grothen\-dieck Festschrift, vol II,
Birkh\"{a}user, Boston (1990), 249--309).

\item{ [Fo00]} J.-M. Fontaine, {\it Arithm\'etique des repr\'esen\-tations
galoisiennes $p$-adiques}, pr\'{e}\-pu\-bli\-cation, Orsay 2000-24. A
para\^{\i}tre dans Ast\'{e}risque.

\item{  [Fo02]} J.-M. Fontaine, {\it Presque-${\Bbb
C}_p$-repr\'{e}\-sen\-ta\-tions}, pr\'{e}publication, Orsay 2002-12.

\item{ [FP02]} J.-M. Fontaine et J\'{e}r\^{o}me Pl\^{u}t, {\it Espaces de
Banach-Colmez}, en pr\'{e}paration.

\item{ [Ke02]}  K. Kedlaya, {\it A $p$-adic local monodromy theorem},
preprint, Berkeley (2001).

\item{ [Me02]}  Z. Mebkhout, {\it Analogue p-adique du
th\'{e}or\`{e}me de Turrittin
et le th\'{e}or\`{e}me de la mo\-no\-dro\-mie p-adique},  Inv. Math.
{\bf 148} (2002),
319--351.

\item{ [Sen80]} S.Sen, Continuous Cohomology and $p$-adic Galois
Representations, Inv. Math. {\bf 62} (1980), 89--116.

\item{ [Ts98]} N. Tsuzuki, {\it Slope filtration of quasi-unipotent
overconvergent $F$-isocrystals}, Ann. Inst. Fourier {\bf 48} (1998),
379--412.

\item{ [Wa96]} N. Wach, {\it Repr\'{e}sentations $p$-adiques
potentiellement cristallines}, Bull. S.M.F. {\bf 124} (1996), 375--400.

\end